\documentclass[12pt,a4paper,makeidx]{amsart}
\usepackage[latin1]{inputenc}        % accents
\usepackage[dvips]{graphics}
\usepackage[dvips]{graphicx}
\usepackage{mathrsfs}
\let\mathcal\mathscr
\usepackage{amssymb}
\usepackage{amsmath}
\usepackage{amsfonts}
\usepackage{latexsym}
\usepackage[T1]{fontenc}
\usepackage[latin1]{inputenc}
\usepackage{hyperref}
\usepackage[all,ps]{xy}
\RequirePackage{amsthm}
\RequirePackage{amssymb}
\usepackage[all,ps]{xy}
\usepackage{latexsym}
\usepackage{amscd}
\usepackage[latin1]{inputenc}
\usepackage[T1]{fontenc}
\usepackage{color}
\RequirePackage{amsthm}
\RequirePackage{amssymb}
\usepackage[all,ps]{xy}
\usepackage{latexsym}

\def\B{{\bf B}}

\def\N{{\bf N}} 
\def\P{{\bf P}}

\def\Q{{\bf Q}}

\def\Betti{\mathop{\rm Betti}\nolimits}

\def\codim{\mathop{\rm codim}\nolimits}

\def\Exc{\hbox{\rm Exc}}

\def\Ii{\mathop{\rm Iitaka}\nolimits}

\def\Im{\mathop{\rm Im}\nolimits}

\def\lra{\longrightarrow}

\def\mult{\mathop{\rm mult}\nolimits}

\def\Supp{\mathop{\rm Supp}\nolimits}
\def\Sym{\mathop{\rm Sym}\nolimits}
\def\Var{\mathop{\rm Var}\nolimits}
\def\vol{\mathop{\rm vol}\nolimits}

\def\dra{\dashrightarrow}
\def\tilde{\widetilde}
\def\eps{\varepsilon}
\def\phi{\varphi}
\def\a{{\alpha}}

\def\d{{\delta}}
\def\e{{\varepsilon}}
\def\i{{\infty}}

\def\m{{\mu}}
\def\n{{\nu}}
\def\D{{\Delta}}

\def\cE{{\mathcal E}}

\def\cI{{\mathcal I}}

\def\cO{{\mathcal O}}

\def\cV{{\mathcal V}}

\def\Nklt{\mathop{\textrm {Non-klt}}\nolimits}
\def\Exc{\mathop{\rm Exc}\nolimits}

\def\base{\mathop{\rm Base}\nolimits}

\def\le{\leqslant}

\newtheorem{thm}{Theorem}[section]

\newtheorem{cor}[thm]{Corollary}
\newtheorem{rmk}[thm]{Remark}

\newtheorem{lem}[thm]{Lemma}

\title[On the uniformity of the Iitaka fibration]{On the uniformity of the Iitaka fibration}
\author[G.~Pacienza]{Gianluca Pacienza}
\date{\today}

\begin{document}
\maketitle

\numberwithin{equation}{section}

%%%%%%%%%%%%%%%%%%%%%%%%%
%
\begin{abstract}
We study pluricanonical systems on smooth projective varieties of positive Kodaira dimension, following the approach of Hacon-McKernan, Takayama and Tsuji succesfully used in the case of varieties of general type. 
We prove a uniformity result for the Iitaka fibration
$X\dra \Ii(X)$ of smooth projective varieties of 
positive Kodaira dimension, provided that $\Ii(X)$ is not uniruled, the variation of the fibration is maximal, and the generic fiber has a good minimal model.
\end{abstract}
%
%%%%%%%%%%%%%%%%%%%%%%%%%%
%
%%%%%%%%%%%%%%%%%%%%%%%%
%
\section{Introduction}
%
%%%%%%%%%%%%%%%%%%%%%%%%
Following  Tsuji 
\cite{tsu1} and \cite{tsu2},  Hacon and McKernan \cite{hm}, and Takayama \cite{taka} have independently given an algebro-geometric proof of the following beautiful result:

\begin{thm}[Hacon-McKernan, Takayama, Tsuji]
\label{thm:0} 
For any positive integer $n$, there exists an integer $m_n$ such that for any smooth complex projective
variety $X$ of general type of dimension $n$, the pluricanonical map 
$$
\phi_{mK_X}: X \dra \P H^0(X,\cO_X(mK_X))^*
$$
is birational onto its image, for all $m\geq m_n$.
\end{thm}
The purpose of this paper is to show that the methods used to prove Theorem \ref{thm:0} allow to obtain  a similar  uniformity result concerning the
pluricanonical maps of  algebraic varieties of arbitrary (positive) Kodaira dimension.
 
Before stating the result we need to recall some facts. Thanks to the work of Iitaka, it is well-known that,
if $\kappa(X)>0$, for large $m$ such that 
$h^0(X,mK_X)\not=0$
the images of the rational maps $\phi_{mK_X}$ 
stabilize i.e. they become birationally 
equivalent to a fiber space 
$$
 \phi_{\infty} : X_{\infty}\lra \Ii(X),
$$
such that 
the restriction of $K_X$ to a very general fiber 
$F$ of $\phi_{\infty}$ has Kodaira dimension $0$ 
and $\dim(\Ii(X))=\kappa(X)$.
This fibration is called the Iitaka fibration of $X$  (see \S \ref{SS:Iitaka} for more details). It is natural to ask (cf. \cite[Conjecture 1.7]{hm}) whether the Iitaka fibration 
of  $X$ 
enjoys a uniformity property as in the case of varieties of general type.  When 
$\kappa(X)=1$ such a result has been
proved in \cite[Theorem 6.1]{fm} with a dependence on the smallest integer $b$ such that 
$h^0(F,bK_F)=1$, and on the Betti number $B_{\dim(E')}$ of a non-singular model $E'$ of the cover $E\to F$ of the general fiber $F$ associated to the unique element of $|bK_F|$ (when $X$ is a 3-fold with $\kappa(X)=1$ this extra dependence may be dropped, see \cite[Corollary 6.2]{fm}). Here we generalize the Fujino-Mori result to arbitrary Kodaira dimension, under extra hypotheses.

\begin{thm}\label{thm:main} 
For any positive integers $n,b,k$, there exists an integer $m(n,b,k)>0$ such that, for any algebraic fiber space 
$f:X\to Y$, with $X$ and $Y$ smooth projective 
varieties,  $\dim(X)=n$, with generic fiber $F$ of $f$
of zero Kodaira dimension, 
and such that:
\begin{enumerate}\label{eq:extra}
\item[(i)] $Y$ is not uniruled;
\item[(ii)] $f$ has maximal variation;
\item[(iii)] the generic fiber $F$ of $f$ has a good minimal model; 
\item[(iv)]  $b$ is the smallest integer such that $h^0 (F, bK_F)\not=0$, and ${\Betti}_{\dim(E')}(E')\leq k$, where $E'$ is a non-singular model of the cover $E\to F$ of the general fiber $F$ associated to the unique element of $|bK_F|$;
\end{enumerate} 
then the pluricanonical map 
$$
 \varphi_{mK_X} : X \dra \P H^0(X,\cO_X(mK_X))^*
$$ is birationally equivalent to $f$, for any $m\geq m(n,b,k)$ 
such that $h^0(X, mK_X)\not=0$.
\end{thm}
Recall that when $F$ is a surface, up to a birational transformation, we may assume that the 12th plurigenus  is non-zero and  the 2nd Betti number is bounded by 22. Therefore,
when $\kappa(X)=n-2$, the integer $m(n,b,k)$ only depends on $n$. The existence of good minimal models is known up to dimension $3$ (cf. \cite{k+}). On the other hand, condition (iii) is automatically satisfied for interesting classes of fibrations, e.g. those for which  $c_1(F)$ is zero (or torsion).

The idea to prove Theorem \ref{thm:main}  is quite natural. By the important result proved in \cite{bdpp}, the hypothesis (i) in Theorem \ref{thm:main}  is equivalent to the pseudo-effectivity of the canonical divisor of $Y$. 
Then, a positivity result due to Kawamata (cf. \cite[Theorem 1.1]{kawapos}, where the hypotheses (ii) and (iii) of Theorem \ref{thm:main} appear), for the (semistable part of the) direct image 
%$f_*\cO_X(\n K_{X/Y})$ 
of the relative pluricanonical sheaf allows to reduce  the problem to the study of effective birationality for multiples of adjoint big divisors 
$K_Y+M$, where $M$ is a big and nef $\Q$-Cartier divisor
such that $\n M$ is integral. The hypothesis (iv) of Theorem \ref{thm:main} is needed to have an effective bound on the denominator 
of the $\Q$-divisor $M$.
Then Theorem \ref{thm:main} is a consequence of the following result, which we prove using the techniques of \cite{hm}, \cite{taka}, and \cite{tsu1}, \cite{tsu2}.
\begin{thm}\label{thm:Mnef} 
For any positive integers $n$ and $\n$, there exists an integer $m_{n,\n}$ such that for any smooth complex projective
variety $Y$  of dimension $n$ with pseudo-effective
canonical divisor, and any big and nef $\Q$-Cartier divisor $M$ on $Y$ such that 
%$K_Y+M$ is big and 
$\n M$ is integral,  the pluriadjoint map 
$$
 \varphi_{m(K_Y+M)} : Y \dra \P H^0(Y,\cO_Y(m (K_Y+M)))^*
$$
is birational onto its image, for all $m\geq m_{n,\n}$ divisible by $\n$.
\end{thm}
As for Theorem \ref{thm:0}, the methods do not lead to an effective constant $m_{n,\n}$.

During the preparation of this article E. Viehweg kindly informed me that he and D.-Q. Zhang were also working on a generalization of the Fujino-Mori result. In their interesting preprint \cite{vz} they study the Iitaka fibration for varieties of Kodaira dimension $2$, and obtain in this case the same uniformity result without the hypotheses (i),(ii) and (iii) appearing in Theorem \ref{thm:main} (and with an effectively computable constant). Their same result, in the case of three-folds, has been obtained independently by Ringler \cite{ringler}.

{\bf Acknowledgements.} I am grateful to  S. Boucksom, F. Campana, L. Caporaso, B. Claudon, O. Debarre, A. Lopez, M. Roth, E. Viehweg and D.-Q. Zhang for useful discussions and/or comments. I wish to thank the Dipartimento di Matematica of the Universit\`a Roma Tre, where this work was done,
for the warm hospitality and the stimulating atmosphere.  My stay in Roma was made possible by an "accueil en d\'el\'egation au CNRS", which is
gratefully acknowledged.

%%%%%%%%%%%%%%%%%%%%%%%%
%
\section{Preliminaries}\label{S:prel}
%
%%%%%%%%%%%%%%%%%%%%%%%%

We recall a number of  basic definitions and results 
that will be freely used in the paper.

%%%%%
%
\subsection{Notation and conventions}
%
%%%%%
We work over the field of complex numbers.
Unless otherwise specified, a divisor will be integral and Cartier.  
If $D$ and $D'$ are $\Q$-divisors on a smooth 
variety $X$ we write $D\sim_\Q D'$, and say that 
$D$ and $D'$ are $\Q$-linearly equivalent, if an integral multiple of $D-D'$ is linearly equivalent to zero. We write $D\equiv D'$ when they are numerically equivalent, that is when they have the same degree
on every curve. The notation $D\leq D'$ means that 
$D'-D$ is effective. If $D=\sum a_i D_i$, 
we denote by $[D]$ the integral divisor 
$\sum [a_i] D_i$, where, as usual, $[a_i]$ is the largest integer which is less than or equal to $a_i$. 
We denote by $\{D\}$ the difference $D-[D]$.
A {\it log-resolution} of a divisor $D\subset X$ is a proper
birational morphism of smooth varieties 
$\mu: X'\to X$, such that the support of $\Exc(\mu)+\mu^*D$ has simple normal crossings. 
The existence of log-resolutions is insured by Hironaka's theorem. 
Given a surjective morphism $f:X\lra Y$ of smooth 
algebraic varieties the {\it relative dualizing sheaf} is the invertible sheaf associated to the divisor $K_{X/Y}:=K_X-f^*K_Y$.
An {\it algebraic fiber space} is a surjective morphism 
$f:X\lra Y$ between smooth projective varieties with connected fibers.

%%%%%
%
\subsection{Volumes, big divisors and base loci}
%
%%%%%
Recall that the volume of a line bundle (see \cite[\S 2.2.C]{l1} for a detailed account on the properties of this invariant)
is the number 
$$\vol_X(D):=\limsup_{m\to +\infty} \frac{h^0(X,mD)}{m^n/n!}$$
It is actually a limit, and we have 
$\vol(mD)=m^{\dim(X)}\vol(D)$. Therefore one can define the volume of a $\Q$-divisor $D$ as $\vol(D):=m^{-\dim(X)}\vol(mD)$, where $m$ is an integer 
such that $mD$ is integral. The volume is invariant
by pull-back via a birational morphism. Moreover we have that $\vol(D)>0$ if and only if $D$ is big, and 
$\vol(D)=D^{\dim(X)}$ for nef divisors.
For a singular variety $Y$, we denote by $\vol(K_Y)$ the volume of the canonical divisor of a desingularization $Y'\to Y$ (which does not depend on the choice of $Y'$). 

 If $V$ is a subvariety of $X$, following \cite{elmnp2} one defines the restricted volume as :
 $$ 
  \vol_{X|V}(A):=\limsup_{m\to+\infty}\frac{h^0(X|V,mA)}{m^d/d!}
 $$
 where 
 $$
  h^0(X|V,mA):=\dim \Im(H^0(X, mA)\to H^0(V,mA_{|V})).$$
Again, it is a limit, and we have that $\vol_{X|V}(mD)=m^{\dim(V)}\vol_{X|V}(D)$ (see \cite[Cor. 2.15 and Lemma 2.2]{elmnp2} ).

We will constantly use Kodaira's lemma : a 
$\Q$-divisor $D$ is big if and only if $D\sim_\Q A+E$, where $A$ 
is a $\Q$-ample divisor and $E$ a $\Q$-effective one.

If $|T|$ is a linear system on $X$, its base locus 
is given by the scheme-theoretic intersection
$$
 \base(|T|) := \bigcap_{L\in|T|} L.
$$
 Recall that given a Cartier divisor $L$ on a variety 
 $X$, its stable base locus (see \cite[pp. 127--128]{l1}) is 
 $$
  \B(L):= \bigcap_{m\geq1} \base(|mL|)
 $$
 and its augmented base locus, which has been defined in \cite{elmnp1}, is 
 $$
  \B_+(L):= \B(mL-H)
 $$
for $m\gg0$ and $H$ ample on $X$ (the latter definition is independent of the choice of $m$ and $H$). One checks that $L$ is ample if, and only if, $\B_+(L)=\emptyset$, and $L$ is big if, and only if, $\B_+(L)\not=X$. In the latter case, $X\setminus \B_+(L)$ is the largest open set on which $L$ is ample.
%%%%%
%
\subsection{Iitaka fibration}\label{SS:Iitaka}
%
%%%%%
We follow \cite[2.1.A and 2.1.C]{l1}.
Let $L$ be a line bundle on a projective variety $X$.
The semigroup $\N(L)$ of $L$ is 
$$
 \N(L):=\{m\geq 0 : h^0(X,mL)\not=0\}.
$$
If $\N(L)$ is not zero, then there exists a natural number $e(L)$, called the exponent of $L$, such 
that all sufficiently large elements in $\N(L)$
are multiples of $e(L)$. If $\kappa(X,L)=\kappa\geq 0$, then $\dim (\phi_{m,L}(X))=\kappa$ for all 
sufficiently large $m\in \N(L)$. Iitaka's result is the following.
\begin{thm}[Iitaka fibrations, see \cite{l1}, Theorem 2.1.33, or \cite{Mo}]
Let $X$ be a normal projective variety and $L$ 
a line bundle on $X$ such that $\kappa(X,L)>0$.
Then for all sufficiently large $k\in\N(L)$ there exists
a commutative diagram 
\begin{equation}\label{eq:diagiitaka}
\xymatrix{ X\ar@{-->}^{\phi_{k,L}}[d]&X_{\i,L}\ar^{\phi_{\i,L}}[d]\ar^{u_{\i}}[l]\\
\Im(\phi_{k,L})&\Ii(X,L)\ar@{-->}^{v_{k,L}}[l].}
\end{equation}
where the horizontal maps are birational. One has 
$\dim(\Ii(X,L))=\kappa(X,L)$. Moreover if we set 
$L_\i=u_{\i}^*L$ and $F$ is the very general fiber of 
$\phi_{\i,L}$, we have
$\kappa(F, L_\i |_F)=0$. 
\end{thm}
We will deal only with the case $L=K_X$, and simply write $\Ii(X):=\Ii(X,K_X)$.
Since the Iitaka fibration is determined only 
up to birational equivalence, and the questions we are interested in are of birational nature, we will often 
tacitly assume that $\Ii(X)$ is smooth, and that we have an algebraic fiber space
$X\lra \Ii(X)$.  Notice that as a consequence of the finite generation of the canonical ring proved in \cite{bchm} we have that, for large $m$, the images 
of the pluricanonical maps $\phi_{mK_X}$ are 
all isomorphic to ${\textrm{Proj}}(\bigoplus_{m\geq 0} H^0 (X,mK_X))$.
 
%%%%%
% 
\subsection{Multiplier ideals}
%
%%%%%

If $D$ is an effective $\Q$-divisor 
on $X$ 
one defines the multiplier ideal as follows : 
$$\cI(X,{D}):=\mu_*\cO_{X'}(K_{X'/X}-[\mu^*D])
$$ 
where $\mu:X'\to X$ is a log-resolution of $(X,D)$.

Notice that if $D$ is a $\Q$-divisor with simple normal crossings, then $\cI(X,D)=\cO_X(-[D])$.
If $D$ is integral we simply have
\begin{equation}\label{eq:int}
 \cI(X,D)=\cO_X(-D).
\end{equation}
Again, we refer the reader to Lazarsfeld's book 
\cite{l2} for  a 
complete treatment of the topic. 
We now recall Nadel's vanishing theorem.

\begin{thm}[{see \cite[Theorem 9.4.8]{l2}}]\label{thm:nadel}
Let $X$ be a  smooth  projective variety. Let $D$ be an effective $\Q$-divisor on $X$, and 
$L$  a divisor  on $X$   such that
 $L-D$ is big and nef. Then, for all $i>0$, we have
$$H^i(X,\cI(X,{D})\otimes \cO_X(K_X+L))=0.$$
\end{thm}

%%%%%
% 
\subsection{Singularities of pairs and Non-klt loci}\label{ss:singpairs}
%
%%%%%

Recall that, in the literature, a pair $(X,D)$ 
is a normal variety together with a $\Q$-Weil 
divisor such that $K_X+D$ is Cartier. 
In this paper the situation is much simpler :
the variety will always be smooth and the divisor 
will be an effective Cartier divisor.
A pair $(X,D)$ is Kawamata log-terminal, klt for short,  (respectively non-klt) at a point $x$, if 
$$
 \cI(X,D)_x = \cO_{X,x}\ \ (\textrm{respectively } 
 \cI(X,D)_x \not= \cO_{X,x}).
$$ 
A pair is klt if it is klt at each point $x\in X$. 
A pair $(X,D)$ is log-canonical, lc for short,  at a point $x$, if 
$$
 \cI(X,(1-\eps)D)_x = \cO_{X,x}\ \ \textrm{for all rational number } 0<\eps<1.
$$ 
A pair is lc if it is lc at each point $x\in X$
(for a survey on singularities of pairs and many related resuts, see \cite{ko}). 
We set 
$$
\Nklt (X,D):=
 \Supp(\cO_X/\cI(X,D) )_{reduced}
$$ 
and call it the Non-klt locus of the pair $(X,D)$.

A simple though extremely useful way of producing example of non-klt pairs
is to consider divisors having high multiplicity at a 
given point, since we have \cite[Proposition 9.3.2]{l2}
\begin{equation}\label{eq:mult}
\mult_x (D)\geq\dim(X) \Rightarrow 
\cI(X,D)_x \not= \cO_{X,x}.
\end{equation}
Notice that if $D=\sum a_i D_i$ is an effective $\Q$-divisor with simple normal crossings, the pair 
$(X,D)$ is klt (respectively lc) if, and only if, $a_i<1$
(resp. $a_i\leq 1$). 
We recall two fundamental results describing the effect of 
small perturbations of $D$ on its Non-klt locus.
\begin{lem}\label{lem:irreduc}
Let $X$ be a smooth projective variety, $x_1$ and $x_2$ two distinct points on $X$, and
$D$ an effective $\Q$-divisor such that $(X,D)$
is lc at $x_1$ and non-klt at $x_2$. Let $V$ be an irreducible component of $\Nklt(X,D)$ passing through $x_1$. Let $B\sim_\Q A+E$ be a big  divisor on $X$, with $A$ $\Q$-ample
 and $E$ $\Q$-effective such that $x_1,x_2\not\in \Supp(E)$. Then there exists an effective divisor $F\sim_\Q B$ and, for any arbitrarily small rational $\delta>0$, there exists a unique rational number 
 $b_{\d}>0$ such that:
\begin{enumerate}
\item $(X,(1-\d)D+b_\d B)$ is lc at $x_1$;
\item $(X,(1-\d)D+b_\d B)$ is non-klt at $x_2$;
\item All the irreducible components of 
$\Nklt (X,(1-\d)D+b_\d B)$ containing $x_1$ are contained in $V$. 
\end{enumerate}
Moreover $\liminf_{\d\lra0}b_\d=0$.
\end{lem}
\begin{proof}
See e.g. \cite[Lemma A.3]{pac}. The reader may also  look at
\cite[Lemma 10.4.8]{l2}
and \cite[Th. 6.9.1]{ko}.
\end{proof}
\begin{lem}\label{lem:dimension}
Let $X$ be a smooth projective variety and
$D$ an effective $\Q$-divisor. Let $V$ be an irreducible component of $\Nklt(X,D)$ of dimension $d$. There exists a dense subset $U$ in the smooth
locus of $V$ and a rational number $\e_0: 0<\e_0<1$ such that, for any $y\in U$, any effective $\Q$-divisor $B$ whose support does not contain $V$ and such that
$$
 \mult_y B|_V >d
$$ 
and any rational  number $\e: 0<\e<\e_0$, the locus
$\Nklt (X,(1-\e)D+B)$ contains $y$.
If moreover $(X,D)$ is lc at  the generic poit of $V$ and $\cI (X,D+B)=\cI (X,D)$ away from $V$,
then $\Nklt (X,(1-\e)D+B)$ is properly contained in $V$ in a neighborhood of any $y\in U$.
\end{lem}
\begin{proof}
See e.g.  \cite[Lemma A.4]{pac}. Again, for similar statements, see  \cite[Lemma 10.4.10]{l2} and \cite[Th. 6.8.1]{ko}.
\end{proof}

%%%%%
%
\section{Positivity results for direct images}
%
%%%%%
In this section we collect results concerning some  positivity properties of the direct image of the relative dualizing sheaf that we will use.

\subsection{The semistable part and a canonical bundle formula}
We recall results contained in \cite[\S 2 and 4]{fm}.
Let $f:X \lra Y$ an algebraic fiber space, whose generic fiber $F$ has Kodaira dimension zero.
Let $b$ be the smallest integer such that the $b$-th
plurigenus $h^0(F,bK_F)$ of $F$ is non-zero. Then there exists a  divisor $L_{X/Y}$ 
on $Y$  (which is unique modulo linear equivalence, and which depends only on the birational equivalence class of $X$ over $Y$) such
that, up to 
birationally modify $X$,  we have
\begin{equation}\label{eq:canonical}
 H^0(Y, ibK_Y + iL_{X/Y})= H^0 (X, ibK_X)
 \end{equation}
 for all $i>0$
(the divisor $L_{X/Y}$ is defined by the double dual
$f_* \cO_X(ibK_{X/Y})^{**}$).
Moreover the divisor $L_{X/Y}$ may be written as 
\begin{equation}\label{eq:decomp}
L_{X/Y}=L^{ss}_{X/Y} + \D
\end{equation}
where  $L^{ss}_{X/Y}$ is a
$\Q$-Cartier divisor, 
called the semistable part or the moduli part (which is compatible with base change),
and $\D$ is an effective $\Q$-divisor called the boundary part. The divisor $L_{X/Y}$ coincides with its moduli part when $f$ is semistable in codimension $1$, and 
\begin{equation}\label{eq:nef}
\textrm {$L^{ss}_{X/Y}$ is nef. }
\end{equation}
The previous results (\ref{eq:canonical}), (\ref{eq:decomp}) and (\ref{eq:nef}) are contained in 
Proposition 2.2, Corollary 2.5 and Theorem 4.5 (iii)
of \cite{fm}. The reader may also look at \cite{kolflips} and \cite[\S 4-5]{Mo}.

For our application it is important to bound the denominators of $L^{ss}_{X/Y}$. 
Let $B$ denote the Betti number $B_{\dim(E')}$ of a non-singular model $E'$ of the cover $E\to F$ of the general fiber $F$ associated to the unique element of $|bK_F|$.  
By \cite[Theorem 3.1]{fm} there exists a positive integer $r=r(B)$  such that 
\begin{equation}\label{eq:betti}
\textrm{$r\cdot L^{ss}_{X/Y}$ is 
an integral divisor. }
\end{equation}
\subsection{Maximal variation and bigness of the semistable part}
Let $f:X\lra Y$ be an algebraic fiber space. Recall that the {\it variation} of $f$ is an integer $\Var(f)$ such that
there exists a fiber space $f':X'\lra Y'$ with $\dim(Y')=\Var(f)$,  a variety $\bar{Y}$, a generically surjective morphism $\varrho : \bar{Y}\lra Y'$ and a generically finite morphism $\pi : \bar{Y}\lra Y$ such that the two fiber spaces 
induced by $\varrho$ and by $\pi$ respectively are birationally equivalent.  The fibration $f$ has {\it maximal variation}
if $\Var(f)=\dim(Y)$. Equivalently,
 $f:X\to Y$ has {maximal variation} if there exists a non-empty open subset $U\subset Y$ such that for any $y_0\in U$ the set $\{y\in U : f^{-1}(y) \sim_{\textrm {birational}} f^{-1}(y_0)\}$ is finite.
As proved by Fujino (\cite [Theorem 3.8]{fu}), we always have 
\begin{equation}\label{eq:fujino}
\kappa(Y,L^{ss}_{X/Y})\leq \Var(f).
\end{equation}
On the other hand, by a result due to Kawamata \cite[Theorem 1.1]{kawapos}, if the generic fiber of $f:X\lra Y$ possesses a good minimal model (i.e. a minimal model whose canonical divisor is semiample), then 
\begin{equation}\label{eq:kawa}
\kappa(Y,L^{ss}_{X/Y})\geq \Var(f).
\end{equation}
In particular, we have
\begin{cor}[Kawamata]\label{cor:kawa} 
Let $f:X\lra Y$ be an algebraic fiber space that has maximal variation and such that 
the generic fiber has a good minimal model. Then 
$L^{ss}_{X/Y}$ is big. 
 \end{cor}
Fujino's inequality (\ref{eq:fujino}) implies that the maximality of the variation is a necessary condition for the bigness of $L^{ss}_{X/Y}$.
 
\subsection{Weak positivity}
Viehweg introduced the notion of weak positivity
for torsion-free coherent sheaf $\cE$ on a projective
variety $V$ : if $V_0$ is the largest open subset on which $\cE$ is locally free, the sheaf $\cE$ is weakly positive if there exists a open dense subset $U$ of $V_0$ such that for any ample divisor $H$
on $V$ and any positive integer $a$, there exists
a positive integer $b$ such that the sheaf 
$(\Sym^{ab}\cE|_{V_0})(bH|_{V_0})$ is generated on $U$ by its global sections on $V_0$ (see \cite{viemod} for a detailed discussion of this notion). 
We will make use of the following positivity result for direct images, due to Campana \cite[Theorem 4.13]{ca}, which improves on previous results obtained by Kawamata \cite{kawab}, Koll\'ar \cite{higherkol} and Viehweg \cite{vie} (see also \cite[Proposition 9.8]{lu}).
\begin{thm}[Campana]\label{thm:campana}
Let $f:V'\to V$ be a morphism with connected fibres between smooth projective varieties. Let $\D$ be an
effective 
$\Q$-divisor on $V'$ whose restriction to the generic
fibre is lc and has simple normal crossings. Then, the sheaf 
$$
 f_* \cO_{V'} (m(K_{V'/V}+\D))
$$
is weakly positive for all positive integer $m$ such that $m\D$ is integral.
\end{thm}
Notice that a locally free sheaf is weakly positive if, and only if, it is pseudo-effective.

%%%%%%%%%%%%%%%%%%%%%%%%
%
\section{Extension of  log-pluricanonical forms}\label{S:ext}
%
%%%%%%%%%%%%%%%%%%%%%%%%
In the course of the proof of Theorem \ref{thm:Mnef}
we will need to lift  (twisted) pluricanonical forms 
on a smooth hypersurface to the ambient variety. 

First, we recall Takayama's  extension result  \cite[Theorem 4.5]{taka}   (cf. \cite[Corollary 3.17]{hm}  for the corresponding result, which is a generalization of a former result of Kawamata's \cite{ka}).

\begin{thm}[Takayama]\label{thm:takaext}
Let $Y$  be a smooth projective variety. Let $H\subset Y$ be a smooth irreducible hypersurface.  
Let $L'\sim_\Q A'+E'$ a big divisor on $Y$  with
 \begin{itemize}
 \item[$\bullet$] $A'$ a nef and big $\Q$-divisor  such that $H\not\subset \B_+(A')$;
\item[$\bullet$]  $E'$ an effective $\Q$-divisor whose support does not contain $H$ and such that the pair
 $(H,E'\vert_H)$ is klt.
\end{itemize}
Then the restriction
$$H^0(Y, m(K_Y +H+L'))\lra H^0(H,m(K_H+L'\vert_H))
$$
is surjective for all integer $m\geq0$.
\end{thm}

The precise statement we need is the following.
 
\begin{cor}\label{cor:ext}
Let $Y$  be a smooth projective variety, $M$ an effective and nef integral Cartier divisor on $Y$. Let $H\subset Y$ be a smooth irreducible hypersurface such that $H\not\subset \Supp(M)$.  
Let $L\sim_\Q A+E$ a big divisor on $Y$  with
 \begin{itemize}
 \item[$\bullet$] $A$ a nef and big $\Q$-divisor  such that $H\not\subset \B_+(A)$;
\item[$\bullet$]  $E$ an effective $\Q$-divisor whose support does not contain $H$ and such that the pair
 $(H,E\vert_H)$ is klt.
\end{itemize}
Then the restriction
$$H^0(X, m(K_Y +M+H+L))\lra H^0(H,m(K_H+M\vert_H+L\vert_H))
$$
is surjective for all integer $m\geq0$.
\end{cor}

For other extension results, all inspired by \cite{siunotgt}, the reader may look
at \cite{cl}, \cite{pa} and \cite{var}.

\begin{proof}[Proof of Corollary \ref{cor:ext}]
We want to apply Theorem \ref{thm:takaext}
to $L'=M+L$. We can write 
$$
M+L= (A+M) + E= \textrm{(big and nef) + effective},
$$
as $M$ is nef. The only thing to check  is that 
$$
H\not\subset \B_+(A+M),
$$ 
assuming $H\not\subset \B_+(A)$. 
But this is immediate, since by the nefness of $M$ we have
$$
 \B_+(A+M)\subset \B_+(A)
$$
and we are done.
\end{proof}

%%%%%%%%%%%%%%%%%%%%%%%%
%
\section{Bounding the restricted volumes from below}\label{S:lowerbound}
%
%%%%%%%%%%%%%%%%%%%%%%%%
It is well-known to specialists that a positive lower bound 
to the restricted volumes of a big divisor $A$ on a variety $X$ allows to construct, along the lines of the Angehrn-Siu proof of the Fujita conjecture, a global section of $K_Y+A$ separating two general points on $X$ (cf. \cite[Proposition 5.3]{taka} and \cite[Theorem 2.20]{elmnp1}). Such a lower bound is the object  of the following result.

\begin{thm}\label{thm:key}
Let $Y$ be a smooth projective variety, $M$ an effective and nef  integral Cartier divisor on $Y$, and 
$V\subset Y$ be an irreducible subvariety not contained in the support of $M$. 
Let $L$ be a big divisor on $Y$ and $L\sim_\Q A+E$ a decomposition such that :
\begin{itemize}
\item[(i)] $A$ is an ample $\Q$-divisor ;
\item[(ii)] $E$ is an effective $\Q$-divisor such that $V$ is an irreducible 
component of $\Nklt(Y,E)$ with $(Y,E)$ lc at the general point of $V$.
\end{itemize}
Then :
$
 \ \ \ \ \ \ \ \ \ \vol_{Y|V} (K_Y+M+L) \geq \vol (K_V+M\vert_V).
$
\end{thm}

The proof of the theorem is a fairly easy consequence of the extension result \ref{cor:ext}
when $\codim(V)=1$. In the general case, it also requires, among other things, the use of Campana's
weak positivity result \ref{thm:campana}.
Using the log-concavity property of the restricted volume,   
established in \cite{elmnp2}, we deduce from Theorem \ref{thm:key} the following consequence which will be the key ingredient in the inductive proof of Theorem \ref{thm:2}.

\begin{cor}\label{cor:key}
Let $Y$ be a smooth projective variety. Let $M'$ be an effective and nef  integral Cartier divisor on $Y$  and 
$V\subset Y$ an irreducible  subvariety not contained in the support of $M'$. 
Let $L$ be a big divisor on $X$ and $L\sim_\Q A+E$ a decomposition such that :
\begin{itemize}
\item[(i)] $A$ is an ample $\Q-$divisor ;
\item[(ii)] $E$ is an effective $\Q-$divisor such that $V$ is an irreducible 
component of $\Nklt(Y,E)$ with $(Y,E)$ lc at the
general point of $V$;
\item[(iii)] $K_Y+L$ is big and $V\not\subset \B_+(K_Y+L)$.
\end{itemize}
Then, for any positive integer $\n$, we have
$$
\vol_{Y|V} (K_Y+\frac{1}{\n}M'+L) \geq \frac{1}{\n^{\dim(V)}}\vol (K_V+\frac{1}{\n}M'\vert_V).
$$
\end{cor}
\begin{proof}
 Write 
 $$
  K_Y+\frac{1}{\n}M'+L=\frac{1}{\n}(K_Y+M'+L)+
  (1-\frac{1}{\n})(K_Y+L).
 $$
 By (iii), thanks to the log-concavity property of the restricted volume proved in \cite[Theorem~A]{elmnp2} we have 
 \begin{eqnarray}\nonumber
 \vol_{Y|V} (K_Y+\frac{1}{\n}M'+L)^{1/d} \geq \frac{1}{\n} \vol_{Y|V} (K_Y+M'+L)^{1/d} + 
 (1-\frac{1}{\n}) \vol_{Y|V} (K_Y+L)^{1/d},
 \end{eqnarray}
 where $d=\dim(V)$. Therefore, by Theorem \ref{thm:key}, we obtain
 \begin{eqnarray}\nonumber
 \vol_{Y|V} (K_Y+\frac{1}{\n}M'+L)^{1/d} \geq
 \frac{1}{\n}\vol (K_V+M'\vert_V)^{1/d}\geq \frac{1}{\n}\vol (K_V+\frac{1}{\n}M'\vert_V)^{1/d}
  \end{eqnarray}
  and the corollary is proved.
\end{proof}
\begin{rmk}{\rm{We will apply Corollary \ref{cor:key} to the base $Y$ of the fibration $f:X\lra Y$ and to (multiples of) a divisor $L=K_Y+ \a L_{X/Y}^{ss}$, where $\a$ will be a certain positive rational number. The  pseudo-effectivity of $K_Y$ that appears among the hypotheses of Theorem \ref{thm:Mnef} is therefore needed here to insure the bigness of the sum $K_Y+L$ that appears in Corollary \ref{cor:key}, hypothesis (iii).}}
\end{rmk}
In the following two subsections we prove Theorem \ref{thm:key}.
\subsection{The case $\codim_{Y}(V)=1$}\label{SS:codim=1}
Hypothesis (ii) simply means 
that $V$
appears with 
multiplicity $1$ in $E$. We may therefore take a modification $\mu:Y'\to Y$  such that
the strict  transform $V'$ of $V$ is smooth,
$\mu^*E=V'+F$ has simple normal crossings  
and moreover 
\begin{equation}\label{eq:notin}
V'\not\subset \Supp(F).
\end{equation} 
Take an integer $m_0>0$ such that
$ m_0(\mu^*A+\{F\})$ has integer coefficients.
By (\ref{eq:notin}) the support of this divisor does not contain $V'$, so we have an inclusion 
\begin{equation}\label{eq:incl}
\xymatrix{ H^0(V',m(K_{V'}+\mu^*M|_{V'}))\ar[d]\ar@^{(->}[d]\\
H^0(V',m(K_{V'}+ \mu^*M|_{V'}+( \mu^*A+\{F\})\vert_{V'}))}
\end{equation}
for any integer $m>0$ divisible by $m_0$.
 %and such that $m\mu^*M$ is integral. 
 Since the pair $(Y', \{F\})$ is klt, applying the extension result \ref{cor:ext} to the divisor $\mu^*A+\{F\}$,
and observing that $\mu^*L-[F]=V'+\mu^*A+\{F\}$, 
we get a surjection
\begin{equation}\label{eq:surj}
\xymatrix{ H^0(Y',m(K_{Y'}+ \mu^*M+\mu^*L -[F]  )) 
\ar@{>>}[d]\\
H^0(V',m(K_{V'}+ \mu^*(M)|_{V'}+(\mu^*A+\{F\})\vert_{V'})).}
\end{equation}

In conclusion we have
\begin{eqnarray*} 
h^0(V,m(K_V+M|_V))&=&h^0(V',mK_{V'}+\mu^*M|_{V'})\\
(\ (\ref{eq:incl})+(\ref{eq:surj})\ )
&\le& 
h^0(Y'\vert V' ,m(K_{Y'}+\mu^*M+\mu^*L  -[F]  ))\\
&\le& 
h^0(Y'\vert V' ,m(K_{Y'}+\mu^*M+\mu^*L    ))\\
&=& 
h^0(Y\vert V ,m(K_Y+M+L)) 
\end{eqnarray*}
so  Theorem \ref{thm:key} is proved in this case.

\subsection{The case $\codim_{Y}(V)\geq2$}\label{SS:codim>1}
We follow here Debarre's presentation \cite[\S 6.2]{D}.
We may assume 
$V$ smooth (see \cite[Lemma 4.6]{taka}). Take a log-resolution $\mu=Y'\to Y$ of $E$, and write
$$
 \mu^* E -K_{Y'/Y}=\sum_F a_FF.
$$
By hypothesis $V$ is an irreducible component 
of $\Nklt(Y,E)$ such that $(Y,E)$ is lc at the general point of $V$. This means that 
\begin{itemize}
\item if $V$ is strictly contained in $\mu(F)$, then $a_F<1$ ;
\item if $V= \mu (F)$, then $a_F\leq1$, with equality 
for at least one $F$. 
\end{itemize}
Thanks to the so-called concentration method
due to Kawamata and Shokurov (see \cite[\S 3-1]{kmm}, and \cite[Lemma 4.8]{taka}) one can further
assume that there exists a unique divisor (which will be denoted by $V'$) among the $F$'s such that 
$\mu(V')=V$. Therefore we have the commutative diagram of smooth varieties :
\begin{equation}\label{eq:inj}
\xymatrix{ V'\ar^{f}[d]\ar[r]\ar@^{(->}[r]&Y'\ar^{\mu}[d]\\
V\ar[r]\ar@^{(->}[r]&Y.}
\end{equation}
We set $G:=\sum_{F\not=V'}a_FF$, and write $[G]$
as a difference of effective divisors without common components $[G]=G_1-G_2$ so that  :

\begin{itemize}
\item $G_2$ is $\mu$-exceptional ;
\item $V\not\subset \mu(\Supp(G_1))$. 
\end{itemize}
Hence, for any integer $m>0$, the sheaf  
$\mu_* \cO_{Y'}(-m[G])$ is an ideal sheaf on $Y$ whose cosupport does not contain $V$, so that
\begin{eqnarray}\label{eq:volinj1}
\ \ H^0(Y, m(K_Y+M+L))&\supset&  
H^0(Y, \mu_* \cO_{Y'}(-m[G])(m(K_Y+M+L)))\\
\nonumber
&\cong& H^0(Y', m(\mu^*(K_Y+M+A+E)-[G]))\\
\nonumber
&\cong& H^0(Y', m(K_{Y'}+\mu^*M+V'+\{G\}+\mu^*A)),
\end{eqnarray}
as soon as the divisors on the right-hand side are integral.
Since the pair $(V',\{G\}|_{V'})$ is klt, we can apply 
the extension result Corollary \ref{cor:ext} to the 
divisor 
$$
 L':= \mu^*L -K_{Y'/Y}-V'-[G]\sim_{\Q}\{G\}+\mu^*A
$$
and to the smooth hypersurface $V'\subset Y'$. 
Hence we get a surjection 
$$
 H^0(Y', m(K_{Y'}+\mu^*M+V'+\{G\}+\mu^*A))
 \twoheadrightarrow H^0(V', m(K_{V'}+\mu^*(M)|_{V'}+L'|_{V'})).
$$
%for all $m$ such that $mM$ is integral.
The last surjection together with the injection (\ref{eq:volinj1}) and, again, the fact that the cosupport of $\mu_* \cO_{Y'}(-m[G])$ does not contain $V$, leads to the 
inclusion :
\begin{equation}\label{eq:volinj2}
H^0(V', m(K_{V'}+\mu^*(M)|_{V'}+L'|_{V'}))
\subset H^0(Y|V, m(K_Y+M+L)).
\end{equation}
On the other hand, thanks to Campana's theorem \ref{thm:campana}, one can show that  for a suitable positive integer $m'$ we have :
$$
 H^0(V', m'(K_{V'/V}+\{G\}|_{V'}+f^*A|_V))\not=0
$$
(see \cite[pages 17-18]{D}).
Hence
%, for any integer $m>0$ such that $mM$ is integral, 
we obtain by multiplication an injection 
$$
 H^0 (V, mm'(K_V+M|_V))\hookrightarrow
 H^0 (V', mm'(K_{V'}+f^*(M|_V) +L'|_{V'})). 
$$
The last inclusion, together with (\ref{eq:volinj2})
and the fact that the restricted volumes are limits
complete the proof of Theorem \ref{thm:key}.
\hfill $\Box$
%%%%%%%%%%%%%%%%%%%%%%%%
%
\section{Point separation for big pluriadjoint systems}\label{S:proofs}
%
%%%%%%%%%%%%%%%%%%%%%%%%
From Theorem \ref{thm:Mnef}
it is easy to deduce
the existence of a uniform positive lower 
bound on the volume of  big adjoint linear systems 
$K_Y+M$ with $M$ nef. 
\begin{cor}\label{cor:1}
For any positive integers $n$ and $\n$, any smooth complex projective
variety $Y$  of dimension $n$ and any big and  nef $\Q$-divisor $M$ with $\n\cdot M$ integral, and such that
$K_Y+M$ is big,
we have :
$$
 \vol(K_Y+M) \geq \frac{1}{(\n\cdot m_{n,\n})^{n}}
$$
where $m_{n,\n}$ is as in Theorem \ref{thm:Mnef}.
\end{cor}
\begin{proof}
Let $m_{n,\n}$ be as in Theorem \ref{thm:Mnef}.
Let $m=\n\cdot m_{n,\n}$.
Let $\mu : Y'\to Y$ be the blow-up along the base locus of $|m(K_Y+M)|$.
Then we can write
$$
 \mu^* m(K_Y+M) = |G| + F, 
$$
where $|G|$ is the base-point-free part, and $F$ is the fixed part.  
In particular $G$ is nef, so $\vol(G)=G^n$. In conclusion we have :
\begin{eqnarray}\label{eq:deg}
\vol (K_Y+M)&=& \frac{\vol (\m^* m(K_Y+M))}{m^n}\geq \frac{1}{m^n}\vol(G)
\\ \nonumber &=&\frac{1}{m^n} G^n =
 \frac{1}{m^n} \deg \phi_{|G|} (Y') \geq \frac{1}{m^n}.
\end{eqnarray}
\end{proof}

On the other hand we will see that a 
sort of converse to Corollary \ref{cor:1} is true. Namely, 
assuming the existence of such a lower bound in dimension $<n$,
we will determine an effective multiple of $K_Y+M$ which is birational. The multiple will still depend on its volume but in a very precise way, sufficient to derive
Theorem \ref{thm:Mnef}.

\begin{thm}\label{thm:2} 
Let $n$ and $\n$ be  positive integers. Suppose there exists a positive constant $v$ such that, for any smooth projective
variety $V$ of dimension $<n$ with pseudo-effective canonical divisor, and any big  
and nef $\Q$-Cartier divisor $N$ on $V$ 
such that
% $K_V+N$ is big and 
$\n N$ is integral,
we have $\vol(K_V+N)\geq v$. Then, there exists two positive constants $a:=a_{n,\n}$ and
$b:=b_{n,\n}$ such that, 
for any smooth projective
variety $Y$ of dimension $n$ with pseudo-effective canonical divisor, and any big and
nef $\Q$-Cartier divisor $M$ on $Y$ such that 
%$K_Y+M$ is big and 
$\n M$ is integral, the rational pluriadjoint map 
$$
 \varphi_{m(K_Y+M)} : Y \dra \P H^0(Y,\cO_Y(m (K_Y+M)))^*
$$
is birational onto its image, for all 
$$
m\geq a + {b\over \vol(K_Y+M)^{1/n}} 
$$
such that $mM$ is integral.
\end{thm}

\subsection{Proof of Theorem \ref{thm:2}}
\label{SS:proof2}
The proof follows the approach
adopted by Angehrn and Siu \cite{as} in their
study of the Fujita conjecture (see also \cite[Theorem 5.9]{ko}), with the variations 
appearing in \cite{hm}, \cite{taka} and \cite{tsu1},\cite{tsu2} 
to make it work for big divisors, and it is based
on the following application of Nadel's vanishing theorem.
 
\begin{lem}\label{lem:pointsep}
Let $Y$ be a smooth projective variety. 
Let $M$ (respectively $E$) be a big and nef 
(resp. a pseudo-effective) $\Q$-divisor. 
Let $x_1,x_2$ be two points outside the support of $E$. Suppose there exists a positive rational number 
$t_0$ such that the divisor
$D_0 \sim t_0 (M+E)$  satisfies the following:

\begin{itemize}
\item[(i)] $x_1, x_2 \in \Nklt (Y,D_0)$;
\item[(ii)] $x_1 \textrm{ is an isolated point in } \Nklt (Y,D_0)$.
\end{itemize}
Then, for all integer $m\geq t_0+1$ such that $(m-1)E+mM$ is integral, there exists a section 
$s\in H^0(Y,K_Y+(m-1)E+mM)$ such that $ s(x_1)\not= 0\textrm { and } s(x_2)= 0$.
\end{lem}
\begin{proof}
Take any integer $m\geq t_0+1$ such that $(m-1)E+mM$ is integral.
Notice that
$$
 D_0+(m-t_0-1)E
$$
is an effective $\Q$-divisor, and that
$$
 (m-1)E+mM - (D_0+(m-t_0-1)E) =(m-t_0)M
$$
is a $\Q$-divisor which is big and nef.
Hence by Nadel's vanishing theorem \ref{thm:nadel} we have:
$$
 H^1 (Y, \cI(Y,{D_0+(m-t_0-1)E})\otimes\cO_Y (K_{Y}+(m-1)E+mM))=0.
$$
Set 
$$
 V_0:= \Nklt (Y,D_0+(m-t_0-1)E)
$$ 
and consider the short exact sequence of $V_0\subset Y$ :
$$
 0\lra \cI(Y,{D_0+(m-t_0-1)E})\lra \cO_Y \lra \cO_{V_0}\lra 0.
$$
Tensoring it with $\cO_Y (K_{Y}+(m-1)E+mM)$, and taking cohomology, we thus get a surjection:
$$
 H^0 (Y, K_{Y}+(m-1)E+mM) \twoheadrightarrow  H^0(V_0, (K_{Y}+(m-1)E+mM)|_{V_0}).
$$
Notice that as the points $x_1, x_2$ lie outside the support of $E$, around them we have 
$$
 \Nklt (Y,D_0)=\Nklt (Y,D_0+(m-t_0-1)E),
$$ 
that is, $V_0$ still contains $x_1,x_2$, the former 
as an isolated point. 
In particular,  there exists a section 
$$
 s\in  H^0 (Y, K_{Y}+(m-1)E+mM)
$$ such that 
$$
 s(x_1)\not= 0\textrm { and } s(x_2)= 0.
$$
\end{proof}
Using (\ref{eq:mult}) it is easy to construct a rational multiple of the big divisor $K_Y+M$ satisfying the condition (i) above. The main problem is that its Non-klt
locus may well have positive dimension at $x_1$.
We will then proceed by descending induction and use the lower bound on the restricted volumes proved in Corollary \ref{cor:key}, in order to cut down the dimension of the Non-klt locus at $x_1$ and end up with  a divisor $D_0\sim t_0(K_Y+M)$, with 
$t_0<a+b/(\vol(K_Y+M))^{1/n}$ and satisfying both hypotheses of Lemma \ref{lem:pointsep}.
In the course of the proof we will invoke the following elementary result.
\begin{lem}\label{lem:gen}
Let $Y$ be a smooth projective variety,
and $M$ an effective $\Q$-divisor such that $K_Y+M$
is big. 
Let $V$ be a subvariety passing through a very 
general point of $Y$ and $\phi:V'\to V$ a desingularization. Then the $\Q$-divisor $K_{V'}+\phi^*M$
is big.
\end{lem}
\begin{proof}
Thanks to the existence of the Hilbert scheme we may assume there exists a smooth family 
$\cV\to B$ and a finite surjective morphism
$\Phi:\cV\to Y$ such that its restriction to the 
general fiber gives $\phi:V'\to V$. Take 
an integer $m>0$ such that $mM$ is integral. Since $\Phi$
is finite, and $\Phi^*|m(K_Y+M)|\subset |m(K_\cV+\Phi^*M)|$, the divisor $K_\cV+\Phi^*M$ is big,
and so is its restriction to the general fiber over $B$.
But the normal bundle of any fiber in the family
is trivial, so by adjunction 
we have $(K_\cV)|_{V'}=K_{V'}$ and we are done.
\end{proof}
\begin{proof}[Proof of Theorem \ref{thm:2}]
The proof follows the Angehrn-Siu approach,  as in the case $M=~0$ that was considered in \cite{hm}, \cite{taka} and in  \cite{tsu1}, \cite{tsu2}. 
We will proceed by descending induction on $d\in\{1,\ldots ,n\}$ to produce an effective $\Q-$divisor 
$D_d\sim t_d (K_Y+M)$ such that :

\begin{enumerate}
 \item $x_1, x_2 \in \Nklt (Y,D_d)$;
 \item $(Y,D_d)$ is lc at a non-empty subset  of $\{x_1,x_2\}$, say at $x_1$;
 \item  $\Nklt (Y,D_d)$ has a unique irreducible component $V_d$ passing through $x_1$ and $\dim V_d\leq d$; 
 \item $t_d <t_{d+1}+v_{d+1}$ with
 $v_{d+1}= \n(d+1)(2/v')^{1/{(d+1)}}(t_{d+1}+2)+\varepsilon$, where $v'\in\{v,\vol(K_Y+M)\}$ and $\eps>0$ may be taken arbitrarily small.
\end{enumerate}

(We set $t_{n}=0$).
 Take two very general points $x_1$ and $x_2$ on 
$Y$. Precisely, they must be outside the augmented base locus of $K_Y+M$, the support of the effective divisor $E$
in the Kodaira decomposition of $K_Y+M\sim_\Q A+E$, the sub-locus of $Y$ covered by the images of $\P^1$, and the union of all the log-subvarieties of the pair $(Y,M)$ which are not of log-general type (i.e. subvarieties $Z$ of $Y$ not contained in $M$ and such that 
$K_{Z'}+ \nu^*M$ is not big, where $\nu:Z'\to Z$ is any desingularization on $Z$).
As in the first step of the Angehrn-Siu proof, 
thanks to the bigness of $K_Y+M$, we can pick an effective $\Q$-divisor $D_{n-1}\sim t_{n-1}(K_Y+M)$ which has multiplicity $>n$
 at both points, as soon as $t_{n-1}\leq n2^{1/n}\vol(K_Y+M)^{-1/n} +\varepsilon$ (with $\eps>0$ arbitrarily small). 
Indeed having multiplicity $\geq r$ at $x_1$ and $x_2$
imposes $2  {{n+r-1}\choose{n}}\sim 2 \frac{r^n}{n!}$ conditions.
On the other hand, for $m\gg0$ such that $[m v_n]M$ is integral, the dimension of  $H^0(Y,[m v_n](K_Y+M))$ grows as $
 \vol(K_Y+M)\frac{[mv_n]^n}{n!}$ so we have
$$
 \vol(K_Y+M)\frac{[mv_n]^n}{n!}> \frac{2m^nn^n}{n!}
$$
as soon as $v_n:= n\big(2/\vol(K_Y+M)\big)^{1/n}+\varepsilon$ (with $\varepsilon>0$
arbitrarily small).
In particular, by (\ref{eq:mult}),  we get that 
$\Nklt (Y, D_{n-1})\ni x_1,x_2$ and, up to multiplying by a positive rational number $\leq1$, we can assume
that $(Y,D_{n-1})$ is lc at one of the two points, say at $x_1$. 
Also, up to performing an arbitrarily small perturbation 
of $D_{n-1}$, thanks to Lemma \ref{lem:irreduc}, 
we may  assume there exists a unique irreducible component of $\Nklt (Y, D_{n-1})$ through $x_1$. The base of the induction is therefore completed.

 For the inductive step we proceed as follows.
 Suppose that we have constructed an effective $\Q$-divisor $D_{d}\sim t_{d}(K_Y+M)$ satisfying 
 conditions (1), (2), (3) and (4) above.  
Suppose for simplicity that $D_{d}$
is nklt at $x_2$ and lc at $x_1$ (the other two possibilities, which are treated in the same way, but 
render the discussion more complicated, are discussed in details
in \cite[\S A.3]{pac}, when $M=0$. The general case can be treated in the same way).
Also, we may assume that $x_1$ is a non-singular point of $\Nklt(Y,D_d)$ (if not, a limiting procedure described in \cite[10.4.C]{l2} allows to conclude). 
 Since
the points lie outside the support of $E$ the same is true for $(Y, D_d+tE)$, where $t:=[t_d]+1-t_d$.
Since $K_Y$ is pseudo-effective,  adding to it any positive multiple of
$K_Y+M$ we still get a big divisor. Therefore we apply Corollary \ref{cor:key} to the divisors
$$
L:=([t_d]+1)(K_Y+M)\sim tA + (D_d+tE) \textrm{ and } 
M':=\n M 
$$ 
and get
$$ 
 \vol_{Y|V_d}(K_Y+M+L)= ([t_d]+2)^d\vol_{Y|V_d}(K_Y+M)\geq   \frac{1}{\n^d} \vol(K_{V_d}+M|_{V_d})
$$
where $V_d$ is the irreducible component of 
$\Nklt(Y, D_d+tE)$ through $x_1$. Since 
$x_1$ is general, $V_d$ cannot be contained in $M$. 
Moreover, always by generality of $x_1$  the divisor
$K_{V_d}+M|_{V_d}$ is big (see Lemma \ref{lem:gen}), and $V$ cannot be uniruled (hence its canonical divisor is pseudo-effective, by \cite{bdpp}).
Then, using the hypothesis,
we have 
$$
 \vol(K_{V_d}+M|_{V_d})\geq v.
$$
Now, we want to add 
to $D_d$ an effective $\Q$-divisor of the form 
$v_d(K_Y+M)$ which has multiplicity $>d$  at 
$x_1$, but
chosen among those restricting to a non-zero
divisor on $V_d$.  
Precisely, for small rational $\d>0$, we add to $(1-\d)D_d$
a divisor $G$
equivalent to 
$$
 (d(2/\vol_{Y|V_d}(K_Y+M))^{1/d} +\varepsilon)(K_Y+M)
$$ 
(which is $\leq (\n d(2/v)^{1/d}(t_d+2) +\varepsilon)(K_Y+M)$
by Corollary \ref{cor:key}),
Using Lemma \ref{lem:dimension}, we get a divisor 
$D_{d-1}\sim t_{d-1}(K_Y+M)$ with 
$$
 t_{d-1}\leq t_d+ \n d(2/v)^{1/d}(t_d+2) +\varepsilon
$$
and such that its $\Nklt$ locus contains $x_1,x_2$.
Moreover $G$ can be chosen such that the new divisor $D_{d-1}$ is 
klt around $x_1$ outside the support 
of $G|_{V_d}$, i.e. 
\begin{eqnarray*}
 \Nklt (Y,D_{d-1}) \textrm { has dimension at } x_1 \textrm { strictly lower than}
\dim(V_d).
\end{eqnarray*}
Again, multiplying $D_{d-1}$ by a rational number $\leq 1$ and applying Lemma \ref{lem:irreduc}, also conditions (2) and (3) are satisfied.
The inductive step is thus proved.

In conclusion, we have obtained an effective $\Q$-divisor $D_0\sim t_0(K_Y+M)$ with 
$$
 t_0< a'_{n,\n}+b'_{n,\n}t_{n-1}\leq a'_{n,\n}+ b'_{n,\n}n(2/\vol(K_Y+M))^{1/n} 
$$
whose $\Nklt$ locus contains $x_1$ as an isolated point, and $x_2$. 
Therefore, by Lemma \ref{lem:pointsep}, 
we deduce the existence of a global section
$$
 s\in  H^0 (Y, K_{Y}+(m-1)K_Y+mM)=H^0 (Y, m(K_Y+M))
$$ 
separating the two points, for all $m>a'_{n,\n}+b'_{n,\n} n(2/\vol(K_Y+M))^{1/n}$ divisible by $\n$.
%To conclude the proof observe that 
%\begin{equation}\nonumber\textrm {the integral divisor $(m-1)K_Y+m M$  is effective}, 
%\end{equation}
%as $K_Y+M$ is big and  $M$ is effective.Therefore we get the inclusion
%$$
 %s\in H^0 (Y, K_{Y}+m (K_Y+M))\subset  H^0 (Y, 2m (K_Y+M))
%$$ 
%and we are done.
\end{proof}
%%%%%
%
\subsection{Proof of Theorem \ref{thm:Mnef}}
\label{SS:proof1}
%
%%%%%
\begin{proof}[Proof of Theorem \ref{thm:Mnef}]
The proof is by induction 
on the dimension of the varieties. 
Theorem  \ref{thm:Mnef} holds for $n=1$. Suppose it holds for $n-1$.
From Corollary \ref{cor:1} we deduce the existence of a positive lower  bound :
$$
 \vol(K_V+N) \geq v_{n-1,\n}
$$
for all pairs $(V,N)$ where $V$ is a smooth projective variety of dimension $\leq n-1$, and $N$ is a big and nef $\Q$-divisor on $V$ such that 
%$K_V+N$ is big and 
$\n N$ is integral.
In particular the hypotheses 
of Theorem \ref{thm:2} are fullfilled. 
Notice that here we use in a crucial way the hypothesis that the denominators of the $N$'s  are bounded. Otherwise the righthand side 
in the inequality (\ref{eq:deg}), which is $1/m^n$,
with $m$ divisible by $\n$,
would go to zero for $\n\lra +\i$.

Consider the pairs  
$(Y,M)$, where $Y$ is $n$-dimensional and 
$M$ is a big and nef $\Q$-divisor such that 
%$K_Y+M$ is big and 
$\n M$ is integral.
For those such that 
volume of $K_Y+M$ is bounded from below, say $\vol(K_Y+M)\geq 1$, then Theorem \ref{thm:2} implies that $|m(K_Y+M)|$ separates points, for all
$$
 m> a+b\geq a + b/\vol(K_Y+M)^{1/n}.
$$
such that $m$ is divisible by $\n$.
For those such that 
$\vol(K_Y+M)< 1$, then {\it a priori} 
the quantity $a + b/\vol(K_Y+M)^{1/n}$ may still be arbitrarily large.
This does not occur : using Theorem \ref{thm:2} and projecting down, 
we have that the variety $Y$ is birational to a subvariety of $\P^{2n+1}$
of degree :
\begin{eqnarray*}
 &\leq&  \Big(a+\frac{b}{\vol(K_Y+M)^{1/n}}\Big)^n\vol(K_Y+M)\\
 &=& \Big(a\vol(K_Y+M)^{1/n}+{b}\Big)^n\leq (a+b)^n.
\end{eqnarray*}
Such varieties are parametrized by an algebraic variety (the Chow variety), so thanks to the Lemma
\ref{lem:noethind} below, the volumes of $K_Y+M$ are also bounded from below by a positive 
constant $c_{n,\n}$ (which is not effective!). 
Hence we may take 
$$
m_{n,\n} := [1+a+{b}/{c_{n,\n}^{1/n}}]
$$  
and we conclude that the pluricanonical system 
$|m(K_Y+M)|$ separates two general points for all $m\geq m_{n,\n}$ divisible by $\n$.
\end{proof}

\begin{lem}[see \cite{taka}, Lemma 6.1]\label{lem:noethind}
Consider the Chow variety 
$${\textrm{Chow}}:=\cup_{d\leq d_n}{\textrm{Chow}}_{n,d_n}(\P^{2n+1}).$$ Let $T=\cup_i T$ be the Zariski closure of 
those points in ${\textrm{Chow}}$.
For any $i$, we have $\inf\{\vol(K_{Y_t}+M_t):t\in T_i{\textrm{ and $K_{Y_t}+M$ is big}}\} >0$.
\end{lem}
\begin{proof}
It is of course sufficient to prove the statement for one $T=T_i$. 
We argue by induction on $\dim(T)$. If the dimension is zero, there is nothing to prove. Suppose $\dim(T)>0$. 
We consider the universal family $U\to T$  and $\tilde U\to U$ a desingularization, together with a line bundle  whose restriction to a smooth fiber is $K_{Y_t}+M_t$
. Let $T^o\subset T$ be the open  subset over which the induced map $p:\tilde U \to T$ is smooth. Let $S\subset T$ be the Zariski dense subset whose points correspond to varieties with $K_{Y_t}+M_t$ big. The $S\cap T^o$ is also dense. 
By construction, for any $s\in S\cap T^o$, the fiber  $\tilde Y_s:=p^{-1}(s)$ is a smooth variety with $K_{Y_t}+M_t$ big. By the upper semicontinuity of the  $h^0$, the same is true for {\it every} fiber over $T^o$, and moreover we have $\inf\{\vol(K_{Y_t}+M_t):t\in  T^o\} >0$.
On the other hand, as for the complement $S\cap (T\setminus T^o)$, we invoke the inductive hypothesis and are done. 
\end{proof}

From Theorem \ref{thm:Mnef} we deduce our main result.

\begin{proof}[Proof of Theorem \ref{thm:main}] 
Fix $n,b$ and $k$. 
Let $r:=r(k)=max\{r(B):B\leq k\}$, where $r(B)$ is the integer appearing in 
(\ref{eq:betti})  and $\n:=r\cdot b$.
Set 
$$
m(n,b,k):=max \{m_{\dim(Y),\n} : 1\leq \dim(Y)\leq n\}$$ 
where $m_{\dim(Y),\n}$ is as in Theorem \ref{thm:Mnef}.
Let $f:X\lra Y$ be an algebraic fiber space verifying the hypotheses of Theorem \ref{thm:main}.
Notice that by (\ref{eq:canonical}) and (\ref{eq:decomp}) we have :
\begin{equation}\label{eq:finalinc} 
H^0 (Y, ibK_Y+i L_{X/Y}^{ss})\subset H^0(Y, ibK_Y + iL_{X/Y})= H^0 (X, ibK_X)
\end{equation}
for all $i>0$ divisible by the integer $r$. 
Take $M:={1\over{b} } L^{ss}_{X/Y}$. By  (\ref{eq:nef}) the $\Q$-divisor $M$ is nef and by (\ref{eq:betti}) the divisor $\n \cdot M$ is integral.
The variety $Y$ is non-uniruled, therefore by \cite{bdpp} its canonical divisor $K_Y$ is pseudo-effective. 
Notice moreover that, since by  (\ref{cor:kawa}) 
the divisor $M$ is big, we have that $K_Y+M$ is big. 
 Then, by Theorem  \ref{thm:Mnef}, we get the birationality of the pluriadjoints maps $\varphi_{m(K_Y+M)}$, for all 
 $m\geq m(n,b,k)$ divisible by $\n$.  The inclusion (\ref{eq:finalinc})
 yields the desired uniformity result for the Iitaka fibration of $X$. 
\end{proof}

%%%%%BIBLIOGRAFIA%%%%%%%%
%

\vskip 30pt

\noindent
{\small Gianluca Pacienza\\
Institut de Recherche Math\'ematique Avanc\'ee\\
Universit\'e L. Pasteur et CNRS\\ 
Rue R. Descartes - 67084 Strasbourg Cedex, France \\
E-mail : {\tt pacienza@math.u-strasbg.fr}}

\end{document}